\documentclass{article}
\usepackage{amssymb}
\usepackage{amsfonts}
\usepackage{color}
\usepackage[bottom]{footmisc}

\newtheorem{theorem}{Theorem}

\newtheorem{condition}[theorem]{Condition}

\input{tcilatex}
\begin{document}

\title{Partial-conjugates and Dimensionality of Posets}
\author{Shaofang Qi\thanks{%
E-mail: sqi@syr.edu. I thank Susan Gensemer and Jerry Kelly for helpful
comments and discussions.}}
\maketitle

\begin{abstract}
The Pareto dominance relation of a preference profile is (the asymmetric
part of) a partial order. For any integer $n$, the problem of the existence
of an $n$-agent preference profile that generates the given Pareto dominance
relation is to investigate the \textit{dimension} of the partial order. We
provide a characterization of a partial order having dimension $n$ in
general.
\end{abstract}

\section{Introduction}

Consider the Pareto dominance relation at a profile of strong preferences
defined on a finite set of objects. If the Pareto relation is observed but
we are ignorant about the preference profile, for an integer $n$, when the
Pareto relation can be generated by an $n$-agent preference profile?\medskip

The following observation allows us to rephrase the question. The existence
of an $n$-agent preference profile implies that the same Pareto dominance
relation can also be generated by an $\left( n+1\right) $-agent preference
profile: assign the additional agent to have the same preference relation as
any one of the existing $n$ agents.\footnote{%
Demuynck (2013) also mentions this observation.} We therefore ask that for
any $n$,\ when the minimum number of individuals whose preference profile
can generate a given Pareto dominance relation is (at most) $n$.\medskip

Sprumont (2001) and Echenique and Ivanov (2011) answer the question for $n=2$%
, from different perspectives. Sprumont imposes a set of \textquotedblleft
regularity\textquotedblright\ conditions and works on a rich continuum of
alternatives, which allows him to utilize a set of simple, and intuitive
basic conditions as (part of) a characterization. Recently, Qi (2013) has
extended Sprumont's basic conditions to a characterization for the finite
case. Echenique and Ivanov (2011) require no specific additional structures
on preferences and focus on the case of a finite set of options; they
convert the question into a graph-coloring problem. To address the analogous
question for $n\geq 3$ is what motivates this work.\medskip 

More generally, the question is equivalent to investigate the \textit{%
dimension} of a \textit{partial order}.\footnote{%
All terminologies will be formally defined in Section 2.} A \textit{partial
order} is a reflexive, antisymmetric, and transitive binary relation defined
on a set of options.\footnote{%
Some authors require irreflxivity in defining partial orders (e.g. Dushnik
and Miller (1941)). But since almost all later work on dimension theory
imposes reflexivity, we follow them and define a partial order to be
reflexive. For an exposition on dimension theory, see the book by Trotter
(1992).} So a Pareto dominance relation plus the diagonal of the binary
relation (i.e. with reflexivity)\ is a partial order. Dushnik and Miller
(1941) introduce the concept of the dimension of a partial order, which is
the minimum number of linear orders whose intersection is the partial order.
The characterizations of 2-dimensional partial orders have been
well-documented. Besides the work mentioned above, there have been other
different characterizations for the 2-dimensional case (see for instance,
Dushnik and Miller (1941), Baker et al. (1972), Kelly (1977), and Trotter
and Moore (1976)). The problem of determining the dimension of a poset
having dimension (at most) $n$ for any fixed $n\geq 3$ is NP-complete
(Yannakakis (1982)).\medskip

The characterization we build on is from Dushnik and Miller (1941) for
2-dimensional partial orders. They introduce the concept of conjugate of a
partial order which is another partial order defined on the same set of
options such that every two distinct options can be comparable by exactly
one of the two partial orders. We extend their concept of conjugates, in two
steps. We first introduce the concept of partial-conjugates which preserves
the properties similar to those hold by conjugates except that the union of
two partial-conjugates partial orders cannot compare all distinct options.
To incorporate this \textquotedblleft completeness\textquotedblright\
property,\ we then introduce a finite sequence of partial orders which have
the partial-conjugates relation and the union of the partial orders of the
sequence has\ every two distinct options comparable. Our main result
provides a characterization, based on our extensions of conjugates, which
generalizes Dushnik and Miller's theorem about conjugates and dimension 2.
Our characterization result is of an \textquotedblleft
existential\textquotedblright\ nature in the sense that we are not providing
an algorithm that can help to determine the dimension of a poset.\medskip 

The rest of the paper is organized as follows. Section 2 discusses notation
and definitions. Since we extend Dushnik and Miller's result, we present
their related concept and theorem in Section 3.\ In Section 4, we introduce
our concept of partial-conjugates along with other concepts, and present our
characterization result. Section 5 concludes with a discussion.

\section{Notation and Definitions}

Let $X$ be a nonempty, finite set. We call $X$ the \textit{ground set}, and
use $|X|$ to denote the number of elements in $X$. Let $\Delta _{X}$ denote
the diagonal of $X\times X$, that is, $\Delta _{X}:=\left\{ \left(
x,x\right) :x\in X\right\} $.\ A \textit{binary relation} $R$ on $X$ is a
nonempty subset of $X\times X$, and we write $xRy$ instead of $\left(
x,y\right) \in R$. A binary relation $R$ on $X$ is \textit{reflexive} if $%
xRx $ for any $x\in X$, \textit{complete} if either $xRy$ or $yRx$ or both
for any $x,y\in X$, \textit{antisymmetric} if $xRy$ and $yRx$ imply that $x,y
$ are identical for any $x,y\in X$, and \textit{transitive} if $xRy$ and $yRz
$ imply $xRz$ for any $x,y,z\in X$.\footnote{%
Note that completeness implies reflexivity. Some authors define completeness
only for any two distinct options.} If $R$ is both reflexive and transitive,
we call it a \textit{quasi-order}. An antisymmetric quasi-order is a \textit{%
partial order}. (That is, a partial order is a reflexive, transitive, and
antisymmetric binary relation.) A complete quasi-order is a \textit{weak
order}. (That is, a weak order is a complete and transitive binary
relation.)\ A complete partial order is a \textit{linear order}. (That is, a
linear order is a complete, transitive, and antisymmetric binary relation.)\
In addition, \textquotedblleft $xRy$ and $yRz$\textquotedblright\ is
shortened to \textquotedblleft $xRyRz$,\textquotedblright\ with a similar
convention applied to any finite conjunctions.\ Let $T_{R}$ denote the 
\textit{transitive closure} of $R$: $xT_{R}y$ if and only if there exist a
positive integer $K$ and elements $x_{1},\dots ,x_{K}$\ such that $%
xRx_{1}Rx_{2}\cdots Rx_{K}=y$. An ordered pair $\left( X,R\right) $ is
called a \textit{partially ordered set}, or simply, a \textit{poset}, if $R$
is a partial order on $X$. Throughout the rest of this paper, a generic
partial order is denoted by $P$. And we use $\mathbb{R}^{n}$ to denote $n$%
-dimensional Euclidean space.\medskip

Let $\left( X,P\right) $ be a poset and consider any elements $x,y\in X$. We
say that $x$ and $y$ are \textit{comparable in }$P$, or simply, \textit{%
comparable}, if either $xPy$ or $yPx$ or both. Accordingly, we say $x$ and $%
y $ are \textit{incomparable in} $P$, or simply, \textit{incomparable}, if $x
$ and $y$ are not comparable in $P$. We write $xNy$ in $P$ if $x$ and $y$
are incomparable in $P$. The \textit{dual} of a partial order $P$ on $X$ is
denoted by $P^{d}$ and is defined by $xP^{d}y$ if and only if $yPx$. The 
\textit{dual} of a poset $\left( X,P\right) $ is denoted by $\left(
X,P^{d}\right) $. Finally, the \textit{dimension} of a poset $\left(
X,P\right) $, denoted $\dim \left( X,P\right) $, is the smallest number of
linear orders (defined on $X$) whose intersection is $P$. It is obvious that
a poset and its dual have the same dimensionality.

\section{Conjugate and Dimension $2$}

Dushnik and Miller (1941) introduce the concept of \textit{conjugate}, which
we illustrate next:

\paragraph{Example 1 (Conjugate).}

Suppose $X=\left\{ x,y,z\right\} $. Consider two partial orders $P$ and $Q$
in Figure 1, both of which are defined on $X$.%
\[
\begin{tabular}{llll}
& $x$ & $y$ & $z$ \\ 
&  &  &  \\ 
$x\medskip $ & $P$ & $P$ &  \\ 
$y\medskip $ &  & $P$ &  \\ 
$z\medskip $ &  &  & $P$ \\ 
\multicolumn{4}{l}{Partial Order $P$}%
\end{tabular}%
\,\,\,\,\,\,\,\,\,\,%
\begin{tabular}{llll}
& $x$ & $y$ & $z$ \\ 
&  &  &  \\ 
$x\medskip $ & $Q$ &  & $Q$ \\ 
$y\medskip $ &  & $Q$ & $Q$ \\ 
$z\medskip $ &  &  & $Q$ \\ 
\multicolumn{4}{l}{Partial Order $Q$}%
\end{tabular}%
\]%
\[
\text{Figure 1: a Partial Order and a Conjugate.} 
\]

$P$ and $Q$ are related in the following sense: (i) if any two distinct
options is comparable in $P$ (resp., $Q$), then it is incomparable in $Q$
(resp., $P$); and (ii) every two distinct options are comparable in either $%
P $ or $Q$. For example, for distinct options $x,y$, $xPy$ but $xNy$ in $Q$.
For distinct options $x,y$; $y,z$; and $x,z$: $xPy$, $yQz$ and $xQz$.\
Additionally, $P\cup Q$ is a linear order on $X$: besides containing the
diagonal $\Delta _{X\times X}$,$\ x\left( P\cup Q\right) y\left( P\cup
Q\right) z$.\medskip

Dushnik and Miller (1941) use \textit{conjugate} to generalize the
relationship of $P$ and $Q$ in Example 1.

\paragraph{Definition (Conjugate, Dushnik and Miller (1941)).}

Let $\left( X,P\right) $ and $\left( X,Q\right) $ be two posets with the
same ground set. $P$ and $Q$ are called \textit{conjugate} partial orders if
every two distinct options of $X$ is ordered in exactly one of them.\bigskip

By definition, for two posets $\left( X,P\right) $ and $\left( X,Q\right) $,
if $P$ and $Q$ are conjugate partial orders, then $P$ and $Q^{d}$ are also
conjugate partial orders, where $Q^{d}$ is the dual of $Q$. The following
lemma generalizes the implication of two conjugate partial orders in Example
1.

\paragraph{Lemma (Lemma 3.51, Dushnik and Miller (1941)).}

Let $\left( X,P\right) $ and $\left( X,Q\right) $ be two posets with the
same ground set $X$. If $P$ and $Q$ are conjugate partial orders, then $%
P\cup Q$ is a linear order defined on $X$.\bigskip

We summarize the properties of partial orders $P$ and $Q$ defined on $X$
that are conjugates:

\begin{condition}
$P$ and $Q$ cannot both order the same two distinct options of $X$.
\end{condition}

\begin{condition}
$P\cup Q$ is a linear order.
\end{condition}

\begin{condition}
$P\cup Q^{d}$ is a linear order.
\end{condition}

Dushnik and Miller provide three characterizations of 2-dimensional partial
orders, one of which connects the dimensionality of 2 to the existence of
conjugate. Our work extends their characterization to $n$-dimensional
partial orders; for comparison, we present their result here.

\paragraph{Theorem (Theorem 3.61 (1) and (3), Dushnik and Miller (1941)).}

Let $\left( X,P\right) $ be a poset. Then $\dim \left( X,P\right) \leq 2$ if
and only if $P$ has a conjugate partial order.

\section{Partial-conjugate and Dimensionality}

We extend the conjugate concept and use the extended concept to characterize 
$n$-dimensional partial orders in general. Our characterization has an
intuition that relates to the natural order defined on a subset of $\mathbb{R%
}^{n}$. We use a poset $\left( X,P\right) $ with $X\subseteq \mathbb{R}^{3}$
to illustrate.

\paragraph{Example 2.}

Let $X=\{(4,2,2),(2,1,4),(1,4,1),(5,3,6),(3,6,5),(6,5,3)\}\subseteq \mathbb{R%
}^{3}$. For convenience, we denote these six elements in $X$ by letters $a$, 
$b$, $c$, $x$, $y$, and $z$:%
\[
\begin{tabular}{|l|l|}
\hline
$a$ & $\left( 4,2,2\right) $ \\ \hline
$b$ & $\left( 2,1,4\right) $ \\ \hline
$c$ & $\left( 1,4,1\right) $ \\ \hline
$x$ & $\left( 5,3,6\right) $ \\ \hline
$y$ & $\left( 3,6,5\right) $ \\ \hline
$z$ & $\left( 6,5,3\right) $ \\ \hline
\end{tabular}%
\]
When we need to specify the ith coordinate of an element a letter denotes,
we use the subscript $i$ for $i\in \left\{ 1,2,3\right\} $. For instance, $%
a=\left( a_{1},a_{2},a_{3}\right) $ where $a_{1}=4$, $a_{2}=2$, and $a_{3}=2$%
. Consider an order $P$ on $X$ such that the diagonal $\Delta _{X}\subseteq
P $ and for distinct options $u,v\in X$, $uPv$ if and only if $u_{i}>v_{i}$
for all $i=1,2,3$, where the symbol $>$ denotes the natural order
\textquotedblleft larger than\textquotedblright\ on $\mathbb{R}$. We
summarize $P$ in Figure 2.%
\[
\begin{tabular}{lllllll}
& $x$ & $y$ & $z$ & $a$ & $b$ & $c$ \\ 
&  &  &  &  &  &  \\ 
$x\medskip $ & $P$ &  &  & $P$ & $P$ &  \\ 
$y\medskip $ &  & $P$ &  &  & $P$ & $P$ \\ 
$z\medskip $ &  &  & $P$ & $P$ &  & $P$ \\ 
$a\medskip $ &  &  &  & $P$ &  &  \\ 
$b\medskip $ &  &  &  &  & $P$ &  \\ 
$c\medskip $ &  &  &  &  &  & $P$%
\end{tabular}%
\]%
\[
\text{Figure 2: a Partial Order }P\text{ on }X\subseteq \mathbb{R}^{3}\text{.%
} 
\]

For the poset $\left( X,P\right) $, $\dim \left( X,P\right) >2$; for a
proof, see for example, Sprumont (2001), Example 1 on page 438. Actually, $%
\dim \left( X,P\right) =3$; one can show this either by finding three linear
orders whose intersection is $P$ or by using Hiraguchi's inequality, $\dim
\left( X,P\right) \leq |X|/2$ for $|X|\geq 4$. Given Dushnik and Miller's
theorem, $P$ doesn't have a conjugate. But consider another partial order $Q$
also defined on $X$, where $\Delta _{X}\subseteq Q$ and for distinct options 
$u,v\in X$, $uQv$ if and only if $u_{i}>v_{i}$ for $i=1,2$, and $u_{i}<v_{i}$
for $i=3$. We present $Q$ in the following Figure 3.\medskip

$P\cup Q$ is also a partial order. In particular, for distinct options $%
u,v\in X$, $u\left( P\cup Q\right) v$ if and only if $u_{i}>v_{i}$ for $%
i=1,2 $. Figure 4 depicts $P\cup Q$, where we use $P$ (instead of $P\cap Q$%
)\ to denote the diagonal.

\[
\begin{tabular}{lllllll}
& $x$ & $y$ & $z$ & $a$ & $b$ & $c$ \\ 
&  &  &  &  &  &  \\ 
$x\medskip $ & $Q$ &  &  &  &  &  \\ 
$y\medskip $ &  & $Q$ &  &  &  &  \\ 
$z\medskip $ & $Q$ &  & $Q$ &  & $Q$ &  \\ 
$a\medskip $ &  &  &  & $Q$ & $Q$ &  \\ 
$b\medskip $ &  &  &  &  & $Q$ &  \\ 
$c\medskip $ &  &  &  &  &  & $Q$%
\end{tabular}%
\]%
\[
\text{Figure 3: a Partial Order }Q\text{ Related to }P\text{ in Figure 2.} 
\]

\[
\begin{tabular}{lllllll}
& $x$ & $y$ & $z$ & $a$ & $b$ & $c$ \\ 
&  &  &  &  &  &  \\ 
$x\medskip $ & $P$ &  &  & $P$ & $P$ &  \\ 
$y\medskip $ &  & $P$ &  &  & $P$ & $P$ \\ 
$z\medskip $ & $%
\color{blue}%
Q$ &  & $P$ & $P$ & $%
\color{blue}%
Q$ & $P$ \\ 
$a\medskip $ &  &  &  & $P$ & $%
\color{blue}%
Q$ &  \\ 
$b\medskip $ &  &  &  &  & $P$ &  \\ 
$c\medskip $ &  &  &  &  &  & $P$%
\end{tabular}%
\]%
\[
\text{Figure 4: The Partial Order }P\cup Q\text{.} 
\]

$P\cup Q$ has a conjugate. We use $R$ to denote a conjugate and depict it,
together with $P$ and $Q$, in Figure 5 (again we use $P$, instead of $P\cap
Q\cap R$,\ to denote the diagonal). $R$ is the partial order such that for
distinct options $u,v\in X$, $uRv$ if and only if $u_{1}>v_{1}$ and $%
u_{2}<v_{2}$. Therefore for distinct $u,v$, $u\left( P\cup Q\cup R\right) v$
if and only if $u_{1}>v_{1}$: $\left( P\cup Q\right) \cup R$ is a linear
order. 
\[
\begin{tabular}{lllllll}
& $x$ & $y$ & $z$ & $a$ & $b$ & $c$ \\ 
&  &  &  &  &  &  \\ 
$x\medskip $ & $P$ & $%
\color{red}%
R$ &  & $P$ & $P$ & $%
\color{red}%
R$ \\ 
$y\medskip $ &  & $P$ &  &  & $P$ & $P$ \\ 
$z\medskip $ & $%
\color{blue}%
Q$ & $%
\color{red}%
R$ & $P$ & $P$ & $%
\color{blue}%
Q$ & $P$ \\ 
$a\medskip $ &  & $%
\color{red}%
R$ &  & $P$ & $%
\color{blue}%
Q$ & $%
\color{red}%
R$ \\ 
$b\medskip $ &  &  &  &  & $P$ & $%
\color{red}%
R$ \\ 
$c\medskip $ &  &  &  &  &  & $P$%
\end{tabular}%
\]%
\[
\text{Figure 5: The Partial Order }P\cup Q\text{.} 
\]

We found that the partial orders $P$ and $Q$ preserve a similar flavor to
the idea \textquotedblleft conjugates.\textquotedblright\ In particular, $P$
and $Q$ don't contain any common two distinct options, that is, condition 1
(in Section 3) of conjugate is satisfied. Although under $P\cup Q$, not all
distinct options are comparable, $P\cup Q$ is a partial order. That is, if
condition 2 of conjugate is extended to \textquotedblleft partial
order,\textquotedblright\ $P$ and $Q$ will satisfy it.\ Finally, $P\cup
Q^{d} $ satisfies a similar but not identical extension: $P\cup Q^{d}$ is
not a linear order, but its transitive closure, $T_{P\cup Q^{d}}$, is a
partial order. We generalize the idea in the following definition.

\paragraph{Definition 1 (Partial-conjugate).}

Let $\left( X,P\right) $ and $\left( X,Q\right) $ be two posets with the
same ground set. $Q$ is called a \textit{partial-conjugate of }$P$ if:%
\newline
(i) every two distinct options of $X$ is ordered in at most one of them;%
\newline
(ii) $P\cup Q$ is a partial order;\newline
(iii) $T_{P\cup Q^{d}}$, the transitive closure of $P\cup Q^{d}$, is a
partial order.

\paragraph{Remark.}

If $Q$ is a partial-conjugate of $P$, then $P$ is also a partial-conjugate
of $Q$.\footnote{%
To see this, note that $\left( T_{Q\cup P^{d}}\right) ^{d}=T_{\left( Q\cup
P^{d}\right) ^{d}}=T_{Q^{d}\cup P}$. Since $T_{Q^{d}\cup P}$ is a partial
order, given that $Q$ is a partial-conjugate of $P$, $T_{Q\cup P^{d}}$ is
also a partial order.}\bigskip

In Definition 1, we list conditions (i), (ii), and (iii) analogous to
conditions 1, 2, and 3 in Section 3. Similar to the conditions in Section 3,
the three conditions here are not independent (condition (ii) and (iii)
together will imply condition (i)). Condition (i) preserves condition 1 of
conjugate (in Section 3) and requires empty intersection of a partial order
and its partial-conjugates on comparing any two distinct options. Condition
(ii) extends condition 2 of conjugate\ in the sense that the union of a
partial order and its partial-conjugate satisfies transitivity but not
necessarily completeness. Similarly, condition (iii) extends condition 3 of
conjugate and requires the union of a partial order and the dual of its
partial-conjugate to be transitive in the weaker sense that the transitive
closure of the union is a partial order. Our next definition completes the
extension of conjugate concept to use a sequence of partial orders having
partial-conjugates relation so that all distinct options can be ordered
under the union of the partial orders of the sequence.

\paragraph{Definition 2 (Sequence of Recursive Partial-conjugates).}

Let $\left( X,P_{1}\right) $, \ldots , $\left( X,P_{n}\right) $ be a
sequence of posets with the same ground set. $P_{1},\dots ,P_{n}$ is called
a \textit{sequence of recursive partial-conjugates} if:\newline
(i) for any $k$ such that $2\leq k\leq n-1$, $P_{k}$ is a partial-conjugate
of $\cup _{i=1}^{k-1}P_{i}$;\newline
(ii) $P_{n}$ is a conjugate of $\cup _{i=1}^{n-1}P_{i}$.\bigskip

For instance, in Example 2, the sequence of three partial orders, $%
P_{1},P_{2},P_{3}$, where $P_{1}=P$, $P_{2}=Q$, and $P_{3}=R$, is a sequence
of recursive partial-conjugates.\medskip

For any poset $\left( X,P\right) $, if $P=P_{1}$ and $P_{1},\dots ,P_{n}$ is
a sequence of recursive partial-conjugates, it is possible to split a
partial order of the sequence, say $P_{2}$, into two partial orders that are
partial-conjugates, and the new sequence is also a sequence of recursive
partial-conjugates. Therefore, we are more interested in a sequence of
recursive partial-conjugates with the smallest number of partial orders. The
following definition serves this purpose.

\paragraph{Definition 3 (an $n$-fold Partial Order).}

Let $\left( X,P\right) $ be a poset. The partial order $P$ is $n$-fold if $n$
is the smallest integer such that there exists a sequence of recursive
partial-conjugates $P_{1},\dots ,P_{n}$ where $P_{1}=P$.

\paragraph{Remark 1.}

Let $\left( X,P\right) $ be a poset. If $P$ is $n$-fold and $P_{1},\dots
,P_{n}$ is a sequence of recursive partial-conjugates where $P_{1}=P$, then $%
P_{1}\cup P_{2}$ is $\left( n-1\right) $-fold.

\paragraph{Remark 2.}

Let $\left( X,P\right) $ be a poset. If $P$ is $n$-fold and $P_{1},\dots
,P_{n}$ is a sequence of recursive partial-conjugates where $P_{1}=P$, then $%
P_{k}\cup P_{k+1}$ is not a partial order for any integer $k$ such that $%
1<k<n$. (Otherwise, take the union of $P_{k}\cup P_{k+1}$ and the number of
sequence can be reduced by $1$, contradiction to that $P$ is $n$%
-fold.)\bigskip

So a $2$-dimensional partial order is $2$-fold. The partial order in Example
2, which is $3$-dimensional, is $3$-fold.

\paragraph{Theorem 1.}

Let $\left( X,P\right) $ be a poset. Then $\dim \left( X,P\right) =n$ if and
only if $P$ is $n$-fold, i.e.,\newline
A. If $\dim \left( X,P\right) =n$, then $P$ is at most $n$-fold;\newline
B. If $P$ is $n$-fold, then $\dim \left( X,P\right) \leq n$.\bigskip

\subsection{Proof of Theorem 1A}

We show: If $\dim \left( X,P\right) =n$, then $P$ is at most $n$-fold.

\paragraph{Proof.}

Consider a poset $\left( X,P\right) $ and suppose that $\dim \left(
X,P\right) =n$. Since $\dim \left( X,P\right) =n$, there exist $n$ linear
orders $L_{1},\dots ,L_{n}$ such that 
\[
P_{1}=P=L_{1}\cap \cdots \cap L_{n}. 
\]%
In what follows, we will only use $P_{1}$ to denote both $P$ and $P_{1}$.%
\newline
We show that $P_{1}$ is at most $n$-fold by constructing a sequence of
recursive partial-conjugates $P_{1},\dots ,P_{n}$.\newline
Define: 
\[
P_{2}:=L_{1}\cap \cdots \cap \left( L_{n}\right) ^{d} 
\]%
\[
P_{3}:=L_{1}\cap \cdots \cap \left( L_{n-1}\right) ^{d} 
\]%
\[
\vdots 
\]%
\[
P_{n}:=L_{1}\cap \left( L_{2}\right) ^{d}. 
\]%
We show that (i) for any $k$ such that $2\leq k\leq n-1$, $P_{k}$ is a
partial-conjugate of $\cup _{i=1}^{k-1}P_{i}$; (ii) $P_{n}$ is a conjugate
of $\cup _{i=1}^{n-1}P_{i}$, and therefore, $P_{1},\dots ,P_{n}$ is a
sequence of recursive partial-conjugates. For any $k$ such that $2\leq k\leq
n-1$, since%
\begin{eqnarray*}
\cup _{i=1}^{k-1}P_{i} &=&P_{1}\cup P_{2}\cup \cdots \cup P_{k-1} \\
&=&\left( L_{1}\cap \cdots \cap L_{n}\right) \cup \left( L_{1}\cap \cdots
\cap \left( L_{n}\right) ^{d}\right) \\
&&\cup \cdots \cup \left( L_{1}\cap \cdots \cap \left( L_{n-k+3}\right)
^{d}\right) \\
&=&L_{1}\cap \cdots \cap L_{n-k+2}
\end{eqnarray*}%
and%
\[
P_{k}=L_{1}\cap \cdots \cap \left( L_{n-k+2}\right) ^{d}, 
\]%
every pair of distinct options of $X$ is ordered in at most one of them and $%
\left( \cup _{i=1}^{k-1}P_{i}\right) \cup P_{k}=L_{1}\cap \cdots \cap
L_{n-k+1}$, which is a partial order. Additionally, since $\cup
_{i=1}^{k-1}P_{i}=L_{1}\cap \cdots \cap L_{n-k+2}\subseteq L_{n-k+2}$, and $%
P_{k}=L_{1}\cap \cdots \cap \left( L_{n-k+2}\right) ^{d}$, which implies $%
\left( P_{k}\right) ^{d}\subseteq L_{n-k+2}$, $\left( \cup
_{i=1}^{k-1}P_{i}\right) \cup \left( P_{k}\right) ^{d}\subseteq L_{n-k+2}$.
Therefore, $T_{\left( \cup _{i=1}^{k-1}P_{i}\right) \cup \left( P_{k}\right)
^{d}}$, the transitive closure of $\left( \cup _{i=1}^{k-1}P_{i}\right) \cup
\left( P_{k}\right) ^{d}$, is a partial order. So, $P_{k}$ is a
partial-conjugate of $\cup _{i=1}^{k-1}P_{i}$. It is also obvious that $%
P_{n} $ is a conjugate of $\cup _{i=1}^{n-1}P_{i}$ since $\cup
_{i=1}^{n-1}P_{i}=L_{1}\cap L_{2}$ and $P_{n}=L_{1}\cap \left( L_{2}\right)
^{-1}$.\newline
So we have constructed a sequence of recursive partial-conjugates $%
P_{1},\dots ,P_{n}$ where $P_{1}=P$. And therefore, $P$ is at most $n$-fold.
\ \ \ \ \ $\square $

\subsection{Proof of Theorem 1B}

We show: If $P$ is $n$-fold, then $\dim \left( X,P\right) \leq n$.

\paragraph{Proof.}

Since $P$ is $n$-fold, consider a sequence of recursive partial-conjugates $%
P_{1},\dots ,P_{n}$ where $P_{1}=P$. We first show that for any $2\leq k\leq
n-1$, if there exist $m$ linear orders such that 
\[
\cup _{i=1}^{k}P_{i}=L_{1}\cap L_{2}\cap \cdots \cap L_{m} 
\]%
then we can find another linear order, denoted as $L_{m+1}$, such that 
\[
\cup _{i=1}^{k-1}P_{i}=L_{1}\cap L_{2}\cap \cdots \cap L_{m}\cap L_{m+1}%
\text{.} 
\]%
\newline
To see this, suppose $\cup _{i=1}^{k}P_{i}=L_{1}\cap L_{2}\cap \cdots \cap
L_{m}$ for linear orders $L_{1},\ldots L_{m}$. Since $P_{1},\dots ,P_{n}$ is
a sequence of recursive partial-conjugates, $P_{k}$ is a partial-conjugate
of $\cup _{i=1}^{k-1}P_{i}$. By condition (iii) of Definition 1, $T_{\left(
\cup _{i=1}^{k-1}P_{i}\right) \cup \left( P_{k}\right) ^{d}}$, the
transitive closure of $\left( \cup _{i=1}^{k-1}P_{i}\right) \cup \left(
P_{k}\right) ^{d}$, is a partial order. Therefore, it can be extended to a
linear order, denoted as $L_{m+1}$. Since 
\[
\cup _{i=1}^{k}P_{i}=L_{1}\cap L_{2}\cap \cdots \cap L_{m} 
\]%
and%
\[
\left( \cup _{i=1}^{k-1}P_{i}\right) \cup \left( P_{k}\right) ^{d}\subseteq
L_{m+1} 
\]%
we have 
\[
\cup _{i=1}^{k-1}P_{i}=L_{1}\cap L_{2}\cap \cdots \cap L_{m}\cap L_{m+1} 
\]%
given that $\left( P_{k}\right) ^{d}$ is the dual of $P_{k}$. So we have
found another linear order $L_{m+1}$ such that $\cup
_{i=1}^{k-1}P_{i}=L_{1}\cap L_{2}\cap \cdots \cap L_{m}\cap L_{m+1}$. Since $%
P_{1},\dots ,P_{n}$ is a sequence of recursive partial-conjugates, $P_{n}$
is a conjugate of $\cup _{i=1}^{n-1}P_{i}$. Therefore, $\cup
_{i=1}^{n-1}P_{i}$ is at most dimension 2 and there exist two linear orders $%
L_{1}$ and $L_{2}$ such that 
\[
\cup _{i=1}^{n-1}P_{i}=L_{1}\cap L_{2}\text{.} 
\]%
Give the result we have just proved, there exists a third linear order $%
L_{3} $, such that%
\[
\cup _{i=1}^{n-2}P_{i}=L_{1}\cap L_{2}\cap L_{3}\text{.} 
\]%
\newline
Repeating the same process, there exists a number of linear orders $%
L_{4},\ldots ,L_{n}$ such that%
\[
\cup _{i=1}^{n-3}P_{i}=L_{1}\cap L_{2}\cap L_{3}\cap L_{4} 
\]%
\[
\vdots 
\]%
\[
P_{1}=P=L_{1}\cap L_{2}\cap \cdots \cap L_{n} 
\]%
so, $\dim \left( X,P\right) \leq n$. \ \ \ \ \ $\square $

\section{Discussion}

Extending the work by Dushnik and Miller, we introduce some concepts related
to their conjugate idea and provide a characterization of a partial order
having dimension $n$ in general. However, as in Dushnik and Miller (1941)
and pointed out by\ Sprumont (2001), our characterization result is of an
\textquotedblleft existential\textquotedblright\ nature so that finding the
objects (a partial-conjugate and a sequence of recursive partial-conjugates
here) stated in our characterization is not necessarily easier than finding
the dimension of the partial order. Since the characterization of an $n$%
-dimensional partial order for any given number of $n$ has been open, the
current work hopes to shed some light on that question. A characterization
that consists of some explicit and simpler conditions which can be easier to
test and applied remains an interesting, though challenging, problem.

\pagebreak

\end{document}